\documentclass[leqno,12pt]{amsart}
\usepackage{comment}
\usepackage[all]{xy}
\usepackage{amsthm}
\usepackage{amssymb} 
\usepackage{amsmath,amscd} 
\usepackage{mathrsfs}
\usepackage{url}
\usepackage{a4wide}

\makeatletter
\newcommand\footnoteref[1]{\protected@xdef\@thefnmark{\ref{#1}}\@footnotemark}
\makeatother

\def\G{\Gamma}

\def\kk{{\mathbf k}}

\def\en{{{\rm en}}}

\def\A_{{A_e^{\en}}}

\def\frks{{\mathfrak{S}}}

\def\op{\operatorname{op}}

\def\mod{\operatorname{mod}}
\def\Mod{\operatorname{Mod}}

\def\proj{\operatorname{proj}}

\def\Hom{\operatorname{Hom}}
\def\End{\operatorname{End}}

\def\Ext{\operatorname{Ext}}

\def\Ker{\operatorname{Ker}}
\def\Im{\operatorname{Im}}

\def\gldim{\operatorname{gldim}}

\newtheorem{lemma}{Lemma}[section]
\newtheorem{proposition}[lemma]{Proposition}
\newtheorem{theorem}[lemma]{Theorem}
\newtheorem{corollary}[lemma]{Corollary}

{

}

\theoremstyle{definition}

\newtheorem{example}[lemma]{Example}
\newtheorem{definition}[lemma]{Definition}

\theoremstyle{remark}

\newtheorem{remark}[lemma]{Remark}

\title{Hochschild cohomology of $q$-Schur algebras 
}

\author{Mayu Tsukamoto
}
\date{\today}

\address{Department of Mathematics, Graduate School of Science, Osaka City University, 3-3-138 Sugimoto, Sumiyoshi-ku, Osaka 558-8585 JAPAN}
\email{m13sa30m19@st.osaka-cu.ac.jp}

\begin{document}

\begin{abstract}
We compute the Hochschild cohomology of any block of $q$-Schur algebras. We focus the even part of this Hochschild cohomology ring. To compute the Hochschild cohomology of $q$-Schur algebras, we prove the following two results: first, we construct two graded algebra surjections between the Hochschild cohomologies of quasi-hereditary algebras because all $q$-Schur algebras over a field are quasi-hereditary. Second, we give the graded algebra isomorphism of Hochschild cohomologies by using a certain derived equivalence.
\end{abstract}

\maketitle
\section{Introduction}
{\it $q$-Schur algebras} were introduced by Dipper and James {\cite{DJ}} in order to study the modular representation theory of finite general linear groups. There exists a surjection from the quantum general linear group onto the $q$-Schur algebra (for example, see {\cite[Theorem 11.3.1]{PW}}). It is known that the $q$-Schur algebras over a field are quasi-hereditary (cf. {\cite[Theorem 11.5.2]{PW}}). 

Let $\kk$ be a splitting field for the algebras we consider, and $q \in \kk \setminus \{0\}$. Let $\mathscr{H}_q (d)$ be the Hecke algebra of the symmetric group $\mathfrak{S}_d$ with parameter $q$ over $\kk$. $\mathscr{H}_q (d)$ is the associative algebra with generators $T_1, \ldots , T_{d-1}$ and satisfies quadratic relations $(T_i + 1)(T_i -q)=0$. Let $\mathscr{S}_q (n, d) := \End_{\mathscr{H}_q (d)}(\displaystyle\bigoplus_{\lambda \in \Lambda(n, d)} x_{\lambda}\mathscr{H}_q (d))$ be the $q$-Schur algebra associated with $\mathscr{H}_q (d)$ and $\Lambda(n, d)$ (see {\cite[Chapter 4]{Mt}} for the details on the above definitions), where $\Lambda(n, d)$ is the set of all sequences of non-negative integers $(\lambda_1, \cdots, \lambda_n)$ such that $\displaystyle\sum_{i=1}^n \lambda_i = d$ and $x_{\lambda} = \displaystyle\sum_{\pi \in \mathfrak{S}_{\lambda}}T_{\pi}$ where $\mathfrak{S}_{\lambda}$ is the Young subgroup corresponding to $\lambda$. 

$\Lambda(n, d)$ is a poset with the dominance ordering $\preceq$. Let $n<d$ and $\xi_{\lambda} :\displaystyle\bigoplus_{\lambda \in \Lambda(d, d)} x_{\lambda}\mathscr{H}_q (d) \twoheadrightarrow x_{\lambda} \mathscr{H}_q (d) \hookrightarrow \displaystyle\bigoplus_{\lambda \in \Lambda(d, d)} x_{\lambda}\mathscr{H}_q (d)$. Then $\xi_{\lambda} \in \mathscr{S}_q(d, d)$ is an idempotent. We set $\xi:= \displaystyle\sum_{\lambda \in \Lambda(n, d)} \xi_{\lambda}$. Then we obtain $\mathscr{S}_q (n, d) \cong \xi \mathscr{S}_q (d, d) \xi$, and we call the above idempotent $\xi$ the Green's idempotent. In the case $q=1$, $\mathscr{S}_1 (n, d)$ is the Schur algebra.  

The theory of cohomology of associative algebras was introduced by Hochschild {\cite{HH}}. {\it Hochschild cohomology} of associative algebras is important in many areas of mathematics, such as ring theory, geometry, representation theory and so on. For example, it was observed that the second Hochschild cohomology group of an associative algebra $A$ controls the deformation theory of $A$ {\cite{G2}}. 
The Hochschild cohomology is a graded algebra with the Yoneda product. One of the most important properties of Hochschild cohomology is its invariance under derived equivalences, proved by Rickard in {\cite[Proposition 2.5] {R}}.

In general, it is difficult to compute the Hochschild cohomology. For several kinds of algebras, the Hochschild cohomologies are calculated. For example, Benson and Erdmann {\cite{BE}} give the Hochschild cohomology of the Hecke algebra of the symmetric group at a root of unity in characteristic zero. 
In this paper, we compute the Hochschild cohomology of $q$-Schur algebras, following the method of {\cite{BE}}, which we give in \S 3.1. Thus their assumption might change to our assumption, see {\cite{MR2761936}} for more detail. 

The structure of this paper is as follows: 

First, we construct the following two graded algebra surjections between the Hochschild cohomology of quasi-hereditary algebras: 
\begin{enumerate}
\item[(i)] If $S$ is quasi-hereditary and $H$ is a heredity ideal in $S$, then there exists a graded algebra surjection from ${\rm HH}^{\ast}(S)$ onto ${\rm HH}^{\ast}(S/H)$;
\item[(ii)] Let $S$ be a quasi-hereditary algebra. We fix a complete set $\{ L(\lambda) \; | \; \lambda \in \Lambda \}$ of simple $S$-modules and a set of orthogonal idempotent $\{e_{\lambda} \; | \; \lambda \in \Lambda \}$ in $S$. For $\pi \subseteq \Lambda$, we put $\epsilon:= \displaystyle\sum_{\lambda \in \pi} e_{\lambda}$. If $\Lambda \setminus \pi$ is a poset ideal, then $\epsilon S \epsilon$ is a quasi-hereditary algebra and there exists a graded algebra surjection from ${\rm HH}^{\ast}(S)$ onto ${\rm HH}^{\ast}(\epsilon S \epsilon)$. 
\end{enumerate}

Second, we compute the even part of the Hochschild cohomology of $q$-Schur algebras by using the above surjections. In particular, the Green's idempotents of the $q$-Schur algebras $\mathscr{S}_q(n, d)$ satisfies the assumption (ii). Therefore we have the following graded algebra surjection: 
\begin{equation*}
{\rm HH}^{\ast}(\mathscr{S}_q (d, d)) \twoheadrightarrow {\rm HH}^{\ast} (\mathscr{S}_q (n, d)).
\end{equation*} 

Third, we construct an explicit bimodule resolution of a certain block $A_e$ of $q$-Schur algebras, and determine the dimensions of Hochschild cohomology groups.
\begin{eqnarray*}
{\rm dim} \;{\rm HH}^i (A_e) = \left\{
\begin{array}{lll}
e, \; & ( i = 0);\\
1, & ( 1 \leq i \leq 2(e-1)); \\
0, & ( 2(e-1) < i).
\end{array}
\right.
\end{eqnarray*}

Finally, we describe the $\kk$-algebra structure of the even part of the Hochschild cohomology ring of $q$-Schur algebras.     

\section{Hochschild cohomology of quasi-hereditary algebras}
In this section, we construct the following two graded algebra surjections between the Hochschild cohomology of quasi-hereditary algebras.

\subsection{Hochschild cohomology}

Let $R$ be an associative algebra over a commutative ring $K$. First, we recall the definition of the Hochschild cohomology of $R$ 
\begin{eqnarray*}
{\rm HH}^i(R):={\rm Ext}^i_{R^{\en}}(R, R),
\end{eqnarray*}
where $R^{\en}$:$=R \otimes_{K} R^{\op}$ (for example, see \cite[Chapeter XI \S 4]{CE}).
This may be expressed in terms of the standard resolution:
\begin{eqnarray}
 \cdots \xrightarrow{d_{3}} R^{\otimes 4} \xrightarrow{d_{2}} R^{\otimes
  3} \xrightarrow{d_{1}} R^{\otimes 2} \xrightarrow{d_{0}} R \rightarrow 0.
\end{eqnarray}
This is an $R^{\en}$-free resolution of $R$, where $d_{0}$ is the multiplication
map and 
\begin{eqnarray*}
 d_{i}(r_{0}\otimes r_{1}\otimes \cdots \otimes
  r_{i+1})=\sum_{n=0}^{i}(-1)^n r_{0} \otimes \cdots \otimes
  r_{n} r_{n+1} \otimes \cdots \otimes r_{i+1}.
\end{eqnarray*}
Applying ${\rm Hom}_{R^{\en}}(-,R)$ to (2-1), we have the following complex: 
\begin{eqnarray*}
0 \rightarrow {\rm Hom}_{R^{\en}}(R^{\otimes 2}, R) &\xrightarrow{{\rm Hom}(d_{1}, R)}& {\rm Hom}_{R^{\en}}(R^{\otimes 3}, R) \\
\xrightarrow{{\rm Hom}(d_{2}, R)} {\rm Hom}_{R^{\en}}(R^{\otimes 4}, R) &\xrightarrow{{\rm Hom}(d_{3}, R)}& \cdots .
\end{eqnarray*}
Thus we have ${\rm HH}^i (R)= {\rm HH}^i(R, R)={\rm Ker Hom}(d_{i+1}, R)/ {\rm Im Hom}(d_{i}, R)$.

We denote by ${\rm HH}^\ast(R)$:$=\displaystyle\bigoplus_{i \geq 0} {\rm HH}^i(R)$ the Hochschild cohomology ring of $R$, where the multiplication is
given by the Yoneda product (cf.\ \cite{B}). 
We denote by $\star$ the Yoneda product in ${\rm HH}^\ast(R)$. Let $\alpha \in {\rm HH}^i(R)$ and $\beta \in {\rm HH}^j(R)$ be the
elements which are represented by $\alpha \in {\rm Ker Hom}(d_{i+1}, R)$
and $\beta \in {\rm Ker Hom}(d_{j+1}, R)$, respectively. Then $\alpha
\star \beta \in {\rm HH}^{i+j}(R)$ is given as follows. There exists the
following commutative diagram of $R^{\en}$-modules
\begin{equation*}
\xymatrix{
\cdots \ar[r] &P_{i+j} \ar[d]_{\sigma_{i}} \ar[r]^{d_{i+j}} &\cdots \ar[r]^{d_{j+2}} &P_{j+1} \ar[d]_{\sigma_{1}} \ar[r]^{d_{j+1}} &P_{j} \ar[d]_{\sigma_{0}} \ar[dr]^{\beta} &\\
\cdots \ar[r] &P_{i} \ar[r]^{d_{i}} &\cdots \ar[r]^{d_{2}} &P_{1} \ar[r]^{d_{1}}
 &P_{0} \ar[r]^{d_{0}} &R \ar[r] &0,
 }
\end{equation*}
where $\sigma_{t} \; (0\leq t \leq i)$ are liftings of $\alpha$. Then we
have $\alpha \star \beta= \alpha \circ \sigma_{i} \in {\rm HH}^{i+j}(R)$. It is known that $\alpha \star \beta$ is independent of choices of
representatives $\alpha$, $\beta$ and liftings $\sigma_{t} \; (0 \leq t \leq i)$. 
Moreover, with this product, Gerstenhaber  proved that ${\rm HH}^\ast(R)$ is a super commutative algebra in {\cite{G}}. That is, for homogeneous elements $\eta \in {\rm HH}^n(R)$ and $\theta \in {\rm HH}^m(R)$, we have $\eta \star \theta = (-1)^{nm} \theta \star \eta$. 

In particular, ${\rm HH}^{\rm ev}(R)$ is a commutative algebra, where ${\rm HH}^{\rm ev}(R):= \displaystyle\bigoplus_{i \geq 0} {\rm HH}^{2i}(R)$. ${\rm HH}^{\rm ev}(R)$ is called the even part of ${\rm HH}^{\ast}(R)$.

\subsection{Quasi-hereditary algebras}

We recall the definition of the quasi-hereditary algebra. This notion was first introduced by Scott {\cite{S}} to study highest weight categories in the representation theory of semisimple complex Lie algebras and algebraic groups. Cline, Parshall and Scott proved many important results in {\cite{CPS}}, see also {\cite{PS}}. In {\cite{DR}}, for a semiprimary ring, Dlab and Ringel gave another definition of quasi-hereditary by using an ideal chain. Let $S$ be a finite dimensional algebra over an algebraically closed field $\kk$. Let $J(S)$ be the Jacobson radical of $S$. We denote by $S\mod$ the category of finitely generated left $S$-modules. We denote by $S\proj$ the category of finitely generated projective left $S$-modules. 
  
\begin{definition}
Let $H$ be a two-sided ideal of $S$.
If $H$ satisfies the following conditions, we call $H$ a {\it heredity ideal} in $S$:  
\begin{enumerate}
\item[(i)] $HH=H$;
\item[(ii)] ${\rm Hom}_S(H,S/H) =0$;
\item[(iii)] $HJ(S)H=0$.
\end{enumerate}
\end{definition}

\begin{definition}[Cline-Parshall-Scott{\cite{CPS}}, Dlab-Ringel{\cite{DR}}]
$S$ is called a {\it quasi-hereditary algebra} if there exists a
chain of ideals 
\begin{equation*}
S=H_{0}>H_{1}> \cdots >H_{n}=0
\end{equation*}
with $H_{i}/H_{i+1}$ heredity ideals in $S/H_{i+1}$, for $0 \leq i <n$. Such a chain of ideals is called a {\it heredity chain} of $S$. 
\end{definition}

We fix a complete set of pairwise non-isomorphic simple $S$-module $\{ L(\lambda) \; | \; \lambda \in \Lambda \}$ and we fix a partial ordering $\leq$ on the index set $\Lambda$. For $\lambda \in \Lambda$, we write $P(\lambda)$ (resp.\ $I(\lambda)$) for the projective cover (resp.\ injective hull) of $L(\lambda)$.

\begin{definition}
For $\lambda \in \Lambda$, there is a unique maximal submodule $K(\lambda)$ of $P(\lambda)$ which satisfies the following condition:
If $[J(S) P(\lambda)/K(\lambda):L(\mu)] \not = 0$, then we have $\lambda > \mu$. We write $\Delta(\lambda) := P(\lambda)/ K(\lambda)$, and we call $\Delta(\lambda)$ the {\it standard module} corresponding to $\lambda \in \Lambda$.

Similarly, for $\lambda \in \Lambda$, we define $\nabla(\lambda)$ as the maximal submodule of $I(\lambda)$ which satisfies the following condition: 
If $[\nabla(\lambda)/L(\lambda): L(\mu)] \not = 0$, then we have $\lambda < \mu.$
We call $\nabla(\lambda)$ the {\it costandard module} corresponding to $\lambda \in \Lambda$.  
\end{definition}

\begin{definition}
Let $M \in S\mod$. If $M$ has a filtration $M=M_0 > M_1 > \cdots > M_i > M_{i+1} > \cdots > M_n=0$ such that $M_i/M_{i+1} \cong \Delta(\lambda)$ for some $\lambda \in \Lambda$ (resp.\ $\nabla(\lambda)$), 
for all $0 \leq i \leq n-1$, then $M$ is called {\it $\Delta$-filtered module} (resp.\ {\it $\nabla$-filtered module}) and $\Delta(\lambda)$ which is isomorphic to $M_i/M_{i+1}$ for some $i$ is called a filtration factor of $M$.
\end{definition}

\begin{remark} [cf.\ Donkin{\cite[A.1 (7)]{D}}]
For a $\Delta$-filtered module $M$, the element $[M]$ in the Grothendieck group $K_0(S)$ of $S\mod$ corresponding to $M$ can be written as 
\begin{equation*}
[M]= \sum_{\lambda \in \Lambda} m_{\lambda} [\Delta(\lambda)] = \sum_{\mu \in \Lambda} (\sum_{\lambda \in \Lambda} m_{\lambda} [\Delta(\lambda): L(\mu)])[L(\mu)].
\end{equation*}
If $[\Delta(\lambda):L(\mu)] \not= 0$, then we have $\mu \leq \lambda$. Thus the coefficients $m_{\lambda}$ are uniquely determined. In other words, the filtration multiplicities do not depend on the choice of the $\Delta$-filtration of $M$. Similarly, we deduce that $\nabla$-filtration multiplicities do not depend on the choice of filtration. Moreover the length of $\Delta$-filtration (resp.\ $\nabla$-filtration) does not depend on the choice of $\Delta$-filtration (resp.\ $\nabla$-filtration). Thus we denote by $fl(M)$ the length of $\Delta$-filtration of $M$ and denote by $(M : \Delta (\lambda))$ (resp.\ $(M : \nabla (\lambda))$) the filtration multiplicity of $\Delta (\lambda)$ (resp.\ $\nabla(\lambda)$).
\end{remark}

\begin{proposition} [cf.\ Donkin{\cite[Appendix Proposition A2.2 (ii)] {D}}] \label{dnkn}
Let $X, Y \in S\mod$. We assume that $X$ is a $\Delta$-filtered module and $Y$ is a $\nabla$-filtered module. Then we have 
\begin{eqnarray*}
\Ext_S^i(X, Y)= \left\{
\begin{array}{ll}
\mathop{\sum}\limits_{\nu \in \Lambda} (X : \Delta(\nu))(Y : \nabla(\nu)),  &i=0;\\
0,  &i \neq 0.
\end{array}
\right.
\end{eqnarray*}
\end{proposition}
 
From the rest of this section, we assume that $S$ is a quasi-hereditary algebra and we fix an ideal $H$ of $S$ which appears in a heredity chain of $S$. We denote by $\overline{S}$ the quotient algebra of $S$ by $H$.
We define $\mathbf{F}:=\overline{S} \otimes_{S} -$ and $\mathbf{F}^{\rm en}:=\overline{S}^{\rm en} \otimes_{S^{\rm en}} -$. Then we have the following lemma:

\begin{lemma} \label{sat}
Let $\lambda, \mu, \nu \in \Lambda$. 
If $\mathbf{F}(L(\lambda ))=L(\lambda )$ and $\mu < \lambda$, then we have $\mathbf{F}(L(\mu ))=L(\mu )$.
If $\mathbf{F}(L(\lambda ))=0$ and $\lambda < \nu$, then we have $\mathbf{F}(L(\nu ))=0$.
\end{lemma}

\begin{lemma} \label{f}
Let $N \in S \mod$. We assume that $N$ is a $\Delta$-filtered module and $\mathbf{F}(\Delta(\lambda))=0$ for each filtration factor $\Delta(\lambda)$ of $N$. Then we have $\mathbf{F}(N)=0.$
\end{lemma}

\begin{proof}
We show the assertion by induction on $fl(N)$.
In the case $fl(N)=1$, it is clear. 

We suppose $fl(N)>1$. Since $N$ is a $\Delta$-filtered module, there exists $\lambda \in \Lambda$ such that 
\begin{equation*}
0 \rightarrow \Delta (\lambda ) \rightarrow N \rightarrow N/\Delta (\lambda ) \rightarrow 0 
\end{equation*}
is a short exact sequence. By this short exact sequence, we have 
\begin{equation*}
\cdots \rightarrow \mathbf{F}(\Delta (\lambda )) \rightarrow \mathbf{F}(N) \rightarrow \mathbf{F}(N/\Delta (\lambda )) \rightarrow 0.
\end{equation*}
Hence the proof is done by the induction hypothesis. 
\end{proof}

\begin{lemma} \label{F}
We fix $\lambda \in \Lambda$ such that $\mathbf{F}(L(\lambda ))=L(\lambda )$. Let $\iota: K(\lambda ) \hookrightarrow P(\lambda )$ be the inclusion map. Then $\mathbf{F}(\iota ): \mathbf{F}(K(\lambda )) \hookrightarrow \mathbf{F}(P(\lambda ))$ is an injection.
\end{lemma}

\begin{proof}
We denote by $\iota |_{HK(\lambda )}$ the restriction of $\iota$ to $HK(\lambda)$.
First we prove $H \Delta(\lambda )=0$ to show that $\iota |_{HK(\lambda )}$ is a surjection. It follows from the definition of standard modules that if $[\Delta (\lambda ):L(\mu )] \neq 0$, then we have $\mu \leq \lambda.$ Thus we deduce from Lemma \ref{sat} that $\mathbf{F}(L(\mu))=L(\mu)$, for each composition factor of $\Delta(\lambda)$. Since we have $H\Delta(\lambda)=0$, it follows that $H P(\lambda )$ is isomorphic to $\iota(H K(\lambda))$. Therefore we have $\mathbf{F}(\iota)$ is an injection. 
\end{proof}

We write $\mathbb{L}_i \mathbf{F}$ for the $i$-th left derived functor of $\mathbf{F}$. In the case $i=1$, we write $\mathbb{L} \mathbf{F}$.

\begin{lemma}\label{lem4}
Let $W$ be a left $\Delta$-filtered module. Then we have $\mathbb{L}_i \mathbf{F}(W)=0$ for any $i>0$. 
\end{lemma}

\begin{proof}
We show the statement by induction on $i$. Firstly, we show this statement in the case $i=1$ by induction on $fl(W)$. 
Now we consider the case $fl(W)=1$. Then there exists $\lambda \in \Lambda$ such that $\Delta(\lambda)$ is isomorphic to $W$. Thus we obtain 
\begin{equation}
0 \rightarrow K(\lambda) \rightarrow P(\lambda) \rightarrow \Delta(\lambda) \rightarrow 0.
\end{equation}
By this short exact sequence (2-2), we have the following exact sequence: 
\begin{equation*}
0 \rightarrow \mathbb{L} \mathbf{F}(\Delta(\lambda)) \rightarrow \mathbf{F}(K(\lambda)) \rightarrow \mathbf{F}(P(\lambda)) \rightarrow \mathbf{F}(\Delta(\lambda)) \rightarrow 0.
\end{equation*}
We assume that $\mathbf{F}(L(\lambda ))=0$. Then we deduce that if $(K(\lambda): \Delta(\mu)) \not=0$, then we have $\mathbf{F}(L(\lambda))=0$. Thus we obtain $\mathbf{F}(\Delta (\mu ))=0$ and we deduce from Lemma \ref{f} that $\mathbf{F}(K(\lambda ))=0$.
For our aim it is sufficient to show that $\mathbf{F}(L(\mu ))=0$ for $\mu \in \Lambda$, where  $\mu$ satisfies  $(K(\lambda ) : \Delta(\mu)) \neq 0$. Therefore we have $\mathbb{L}\mathbf{F}(\Delta(\lambda))=0$. 
If $\mathbf{F}(L(\lambda))=L(\lambda)$, then it follows from Lemma \ref{F} that $\mathbf{F}(\iota )$ is injective. Thus we deduce $\mathbb{L} \mathbf{F}(\Delta(\lambda))=0$. 

We assume that $fl(W)>1$. Then there exists a standard module $\Delta(\lambda)$ and a factor module $Q$ of $W$ such that
\begin{equation}
0 \rightarrow  \Delta(\lambda) \rightarrow W \rightarrow Q \rightarrow 0
\end{equation}
is a short exact sequence. Since $\mathbb{L} \mathbf{F}(\Delta(\lambda))=0$, we obtain 
\begin{equation*}
0 \rightarrow \mathbb{L} \mathbf{F}(W) \rightarrow \mathbb{L} \mathbf{F}(Q) \rightarrow \mathbf{F}(\Delta(\lambda)) \rightarrow \mathbf{F}(W) \rightarrow \mathbf{F}(Q) \rightarrow 0
\end{equation*}
from the short exact sequence (2-3). Therefore we deduce from the induction hypothesis that $\mathbb{L}\mathbf{F}(M)=0$. 

Secondly, we also show the assertion in the case $i>1$ by induction on $fl(W)$. 
If $fl(W)=1$, then there exists $\lambda \in \Lambda$ such that $\Delta(\lambda)$ is isomorphic to $W$.  
Hence we have $\mathbb{L}_{i+1} \mathbf{F}(\Delta(\lambda)) \cong \mathbb{L}_i \mathbf{F}(K(\lambda))$ from the short exact sequence (2-2). Therefore the assertion follows from the induction hypothesis.

If $fl(W)>1$, then we obtain
\begin{equation*}
\cdots \rightarrow \mathbb{L}_{i+1} \mathbf{F}(\Delta(\lambda)) \rightarrow \mathbb{L}_{i} \mathbf{F}(W) \rightarrow 
\mathbb{L}_i \mathbf{F}(Q) \rightarrow \mathbb{L}_i \mathbf{F}(\Delta(\lambda)) \rightarrow \cdots.
\end{equation*}
Thus we have $\mathbb{L}_{i} \mathbf{F}(W) \cong \mathbb{L}_i \mathbf{F}(Q)$ for any $i \geq 2$. Therefore the assertion follows from the induction hypothesis.
\end{proof}

\begin{lemma} \label{gr}
Let $Y \in \overline{S} \mod$. Then we have $\mathbb{L}_i \mathbf{F}(Y)=0$ for any $i>0$. 
\end{lemma}

\begin{proof}
We show the assertion by induction on $i$. We assume that $\lambda \in \Lambda$ satisfies that $[Y:L(\lambda )] \neq 0$.
For our aim, it is sufficient to show that $\mathbb{L}_i \mathbf{F}(L(\lambda ))=0$ for any $i>0$. If $i=1$, then there exists a submodule $K$ of $\Delta(\lambda )$ such that 
\begin{equation}
0 \rightarrow K \rightarrow \Delta(\lambda ) \rightarrow L(\lambda ) \rightarrow 0
\end{equation}
is a short exact sequence. This short exact sequence (2-4) induces the long exact sequence as follows:
\begin{equation*}
\cdots \rightarrow \mathbb{L} \mathbf{F}(\Delta (\lambda )) \rightarrow \mathbb{L} \mathbf{F}(L(\lambda )) \rightarrow \mathbf{F}(K) \rightarrow \mathbf{F}(\Delta (\lambda )) \rightarrow 0.
\end{equation*} 
Thus we deduce from Lemma \ref{F} that $\mathbb{L} \mathbf{F}(\Delta (\lambda ))=0$. 
If $\mathbf{F}(L(\lambda ))=0$, then $\mathbf{F}(K)=0$. Hence we have $\mathbb{L} \mathbf{F}(L(\lambda ))=0$.
Let $\mathbf{F}(L(\lambda ))=L(\lambda )$. Then we deduce that if $L(\mu)$ is a composition factor of $K$, then we have $\mu < \lambda.$
Thus we obtain from Lemma \ref{F} that  $\mathbf{F}(L(\mu ))=L(\mu )$. It follows from Lemma \ref{sat} that $\mathbf{F}(K) = K$ and $\mathbf{F}(\Delta (\lambda )) = \Delta(\lambda )$. 
Therefore  we deduce $\mathbb{L} \mathbf{F}(L (\lambda ))=0$. 

If $i>1$, then the short exact sequence (2-4) induces the long exact sequence as follows: 
\begin{equation*}
\cdots \rightarrow \mathbb{L}_{i+1} \mathbf{F} (\Delta (\lambda )) \rightarrow \mathbb{L}_{i+1} \mathbf{F}(L(\lambda )) \rightarrow 
\mathbb{L}_i \mathbf{F}(K) \rightarrow \mathbb{L}_i \mathbf{F} (\Delta (\lambda )) \rightarrow \cdots .
\end{equation*}
Then it follows from Lemma \ref{F} that 
\begin{equation*}
\mathbb{L}_{i+1} \mathbf{F}(L(\lambda )) \cong \mathbb{L}_i \mathbf{F}(K).
\end{equation*}
Hence we obtain $\mathbb{L}_{i+1} \mathbf{F}(\Delta (\lambda ))=0$ by the induction hypothesis.
\end{proof}

\begin{lemma}\label{lem2}
Let $X, Y \in \overline{S}\mod$. Then we have  
\begin{equation*}
{\rm Tor}_i^S(X,Y) \cong {\rm Tor}_i^{\overline{S}}(X,Y) 
\end{equation*}
for any $i \geq 0$. 
\end{lemma}

\begin{proof}
For $X \in  \mod \overline{S}$, we define $\mathbf{G}:=X \otimes_{\overline{S}} -$. Then $\mathbf{F}$ and $\mathbf{G}$ are right exact functors, and $\mathbf{F}(P)$ is $\mathbf{G}$-acyclic for any $P \in S \proj$. Hence there exists a Grothendieck spectral sequence (for example, see \cite{CE}) $E = (E^r_{p, q}, E_n)$ of $S\mod$ such that for each $Y \in S\mod$, the following holds: 
\begin{equation*}
E_{p,q}^2 \cong {\rm Tor}^{\overline{S}}_{p}(X,\mathbb{L}_q \mathbf{F}(Y)), E_n = {\rm Tor}_{p+q}^S(X,Y). 
\end{equation*}
Moreover we deduce from Lemma \ref{gr} that $\forall i>0, \mathbb{L}_i \mathbf{F}(Y)=0$ for $i>0$. Hence the assertion holds.  
\end{proof}

\begin{lemma}[cf.\ Cartan-Eilenberg {\cite[Chapter IX Theorem 2.8 (a)] {CE}}] \label{CE} 
Let $\Lambda, \Gamma, \Sigma$ be algebras over $\kk$. 
For $A \in \Mod \Lambda \otimes \Sigma$, $B \in \Lambda \Mod \Sigma, C \in \Gamma \otimes \Sigma \Mod$, we assume that 
\begin{equation*}
{\rm Tor}_n^{\Lambda}(A,B) = 0 = {\rm Tor}_n^{\Sigma}(B,C)
\end{equation*}
for any $n>0$. Then we obtain 
\begin{equation*}
{\rm Tor}_i^{\Lambda \otimes \Sigma}(A \otimes_{\Lambda} B, C) \cong {\rm Tor}_i^{\Lambda \otimes \Gamma}(A, B \otimes_{\Sigma}C)
\end{equation*} 
for any $i \geq 0$. 
\end{lemma}

\subsection{Main result 1}

\begin{theorem}\label{thm2}
Let $S$ be a quasi-hereditary algebra over $\kk$. We assume that $H$ appears in a heredity chain of $S$. 
Then there exists a surjective graded algebra homomorphism:  
\begin{equation*}
\phi :{\rm HH}^{\ast}(S) \twoheadrightarrow {\rm HH}^{\ast}(S/H).
\end{equation*}
\end{theorem}

\begin{proof}
We deduce from Lemma \ref{lem2} and Lemma \ref{CE} that 
\begin{eqnarray*}
{\rm Tor}_i^{S^{\rm en}} (\overline{S}^{\rm en}, S) &\cong & {\rm Tor}_i^S (\overline{S}, \overline{S}^{\op} \otimes_{S^{\op}} S)\\
& \cong & {\rm Tor}_i^{\overline{S}} (\overline{S}, \overline{S}^{\op} \otimes_{S^{\op}} S).
\end{eqnarray*}
Hence we obtain 
\begin{equation*}
{\rm Tor}_i^{S^{\rm en}} (\overline{S}^{\rm en}, S)=0
\end{equation*}
for any $i \geq 1$. Thus we can define the following map for each $i$.
\begin{eqnarray*}
\phi_i : {\rm HH}^i(S) & \twoheadrightarrow & {\rm HH}^i(\overline{S})\\
\hspace{3em} [\alpha] & \mapsto & [\mathbf{F}^{\rm en}(\alpha )]. 
\end{eqnarray*}
Firstly, we show that $\phi_i$ is well-defined.
Let ($P_{\bullet}$, $d_{\bullet}$) be a projective resolution of $S$ as left $S^{\rm en}$-modules. 
We deduce from ${\rm HH}^i(S) \cong \Ker \Hom (d_{i+1}, S)/ \Im \Hom (d_i, S)$ that if $\alpha \in \Im \Hom (d_i, S)$, then there exists $\beta \in \Hom_{S^{\rm en}}(P_i, S)$ such that  $\alpha=\beta \circ d_i$ and the following holds: 
\begin{equation*}
\mathbf{F}^{\rm en}(\alpha )=\mathbf{F}^{\rm en}(\beta \circ d_i)=\mathbf{F}^{\rm en}(\alpha )\circ \mathbf{F}^{\rm en}(d_i).
\end{equation*}
Thus $\mathbf{F}^{\rm en}(\alpha ) \in \Im \Hom(\mathbf{F}^{\rm en}(d_i), \overline{S})$.
Using the fact that $\mathbf{F}^{\rm en}$ is a functor and the Yoneda product is functorial, it follows that $\phi$ is a graded algebra homomorphism. 
It is straightforward to check that $\phi$ is a surjection. 
\end{proof}

\begin{remark}
Let $S$ be a quasi-hereditary algebra with a heredity chain $S=H_{0}>H_{1}> \cdots >H_{n}=0$. Then we deduce the surjective graded algebra homomorphism from ${\rm HH}^{\ast}(S)$ to ${\rm HH}^{\ast}(S/H_{n-1})$ in a different way. Again we take a heredity ideal $H_{n-2}/H_{n-1}$ of $S/H_{n-1}$. Then we have the surjective graded algebra homomorphism from ${\rm HH}^{\ast}(S/H_{n-1})$ to ${\rm HH}^{\ast}(S/H_{n-2})$. Thus we can show Theorem \ref{thm2} by repeating this process inductively.
\end{remark}

We provide an example of a finite dimensional algebra which dose not hold Theorem \ref{thm2}. 

\begin{example}
Let $A$ be the algebra over a field defined by the following quiver
\[
\xymatrix{1\ar@(dl,ul)^{\gamma_{1}} \ar@/^/[r]^{\alpha} &2\ar@(dr,ur)_{\gamma_2}\ar@/^/[l]^{\beta}}
\]
with relations $\gamma_i ^2 \; (i=1, 2), \gamma_2 \alpha, \gamma_1 \beta, \alpha \beta \alpha, \beta \gamma_{2}$. Let $I:= A e_2 A$. Then there dose not exist a surjective graded algebra homomorphism form ${\rm HH}^{\ast}(A)$ to ${\rm HH}^{\ast}(A/I)$. 
\end{example}

We deduce the corollary from Theorem \ref{thm2} as follows:

\begin{corollary} \label{cor1}
Let $S$ be a quasi-hereditary algebra. Then there exists an idempotent $\xi $ of $S$ such that $\xi S \xi$ is a quasi-hereditary algebra and  
 \begin{equation*}
\psi : {\rm HH}^{\ast}(S) \twoheadrightarrow {\rm HH}^{\ast}(\xi S \xi)
\end{equation*}
is a graded algebra surjection.
\end{corollary} 

\begin{proof}
Let $(S_{\xi})^{'}$ be the Ringel dual of $S_{\xi}$. Then there exists a two sided ideal $H^{'}$ of the Ringel dual $S^{'}$ of $S$ such that $H^{'}$ appears in a heredity chain of $S^{'}$, and the following isomorphism is given by the property of the Ringel duality, this is proved by Ringel  in {\cite[Appendix]{Rin}}.  
\begin{equation*}
(\xi S \xi)' \cong S'/H'.
\end{equation*}
Hence we deduce 
\begin{eqnarray*}
{\rm HH}^{\ast}(\xi S \xi ) &\cong& {\rm HH}^{\ast}(S'/H'), \\
{\rm HH}^{\ast}(S) &\cong& {\rm HH}^{\ast}(S').
\end{eqnarray*}  
Moreover we have 
\begin{equation*}
{\rm HH}^{\ast}(S') \twoheadrightarrow {\rm HH}^{\ast}(S'/H') 
\end{equation*}
is a graded algebra surjection from Theorem \ref{thm2}.
Therefore we can construct $\psi$.
\end{proof}

\section{Hochschild cohomology of the $q$-Schur algebras}
In this section we compute the Hochschild cohomology of the $q$-Schur algebras. From the rest of this paper, we denote by $\kk$ a field of characteristic $l \geq 0$. We put 
\begin{equation*}
e:= {\rm inf} \{i \in \mathbb{Z}_{ \geq 1}\; | \; 1+ q + \cdots + q^{i-1} = 0 \; {\rm in}\;  \kk\}.
\end{equation*}
$e$ is called the quantum characteristic. 

\subsection{Preliminaries}
We prepare a certain derived equivalence and an explicit bimodule projective resolution to compute the Hochschild cohomology of $q$-Schur algebras. 
In this subsection, we use some combinatorial notion (e.g.\ $e$-weight, $e$-core and $e$-abacus) to describe the above derived equivalence, so see {\cite[Chapter 2]{W}} for more detail.
Let $B$ a block of $\mathscr{S}_q(n, n)$ of $e$-weight\footnote{\label{fn_label1}{see, for example {\cite[Chapter 2]{W}}}} $w \in \mathbb{Z}_{\geq 0}$ and $B^{'}$ be a block of some $\mathscr{S}_q(m, m)$ with the same $e$-weight $w$. Then $B$ and $B^{'}$ are derived equivalent \cite{CR}. 

\begin{theorem}[{\it Nakayama Conjecture}, Dipper-James {\cite[Theorem 6.7]{DJ}}] \label{block}
The blocks of $q$-Schur algebras $\mathscr{S}_q(n, n)$ are in one-one correspondence with pairs $(w, \tau)$, where $w \in \mathbb{Z}_{\geq 0}$ is an $e$-weight, $\tau$ is an $e$-core\footnoteref{fn_label1} of size $n - we$.
\end{theorem}

From Theorem \ref{block}, we denote by $\mathbf{B}_{\tau, w}$ the block of the $q$-Schur algebra $\mathscr{S}_q (n, n)$ corresponding to the pair $(w, \tau)$.

\begin{theorem} [Chuang-Rouquier {\cite[\S 7.6]{CR}}] \label{CR}
Let $\tau, \tau^{'}$ be $e$-cores. Then the following derived equivalence holds: 
\begin{equation*}
\mathcal{D}^b(\mathbf{B}_{\tau, w} \mod) \simeq \mathcal{D}^b(\mathbf{B}_{\tau^{'}, w} \mod). 
\end{equation*}
\end{theorem}
\begin{definition}
We suppose that $p,w \in \mathbb{Z}_{\geq 0}$ are fixed. A $p$-core $\rho$ is said to be a {\it Rouquier $p$-core} if it has a $p$-abacus\footnoteref{fn_label1} presentation, on which there are at least $w-1$ more beads on runner $i$, than on runner $i-1$,
for $i =1, \ldots ,p-1$. 
\end{definition}
\begin{definition}
Let $w \in \mathbb{Z}_{\geq 0}$ and $\rho$ be a Rouquier core of $e$-weight $w$.
We say that $\mathbf{B}_{\rho, w}$ is a {\it Rouquier block} of a $q$-Schur algebra. 
\end{definition}

The notion of the Rouquier block give important information to us. For example, the following theorem holds. 

\begin{theorem} [Chuang-Miyachi {\cite[Theorem 18]{MR2761936}}] \label{CM} 
Let $\mathbb{F}_l$ be a finite field of $l$ elements. We assume that $l=0$ or $w < l, q \in \mathbb{F}_l $. 
Then $\mathbf{B}_{\rho, w}$ and $B_0 (\mathscr{S}_q (e, e))^{\otimes w} \rtimes \kk \frks_w$ are Morita equivalent, where $B_0 (\mathscr{S}_q (e, e))$ is the principal block of the $q$-Schur algebra $\mathscr{S}_q(e,e)$. 
\end{theorem}

\begin{lemma}[cf.\ {\cite{EN}}]  \label{mori}
We use the same notation in Theorem \ref{CM}.
$B_0 (\mathscr{S}_q (e, e)) $ is Morita equivalent to $A_e $, where $A_e$ is the algebra over a field defined by the following quiver 
\begin{equation*}
Q :=(1) \overset{\alpha(1)}{\underset{\alpha^-(1)}{\rightleftarrows}} \cdots (i-1) \overset{\alpha(i-1)}{\underset{\alpha^-(i-1)}{\rightleftarrows}} (i) \overset{\alpha(i)}{\underset{\alpha^-(i)}{\rightleftarrows}} (i+1) \cdots
\overset{\alpha(e-1)}{\underset{\alpha^-(e-1)}{\rightleftarrows}} (e),
\end{equation*}
with relations 
\begin{eqnarray*}
&\alpha(i)\alpha(i-1), \; \alpha^-(i-1)\alpha^-(i) ,\\
&\alpha(i-1)\alpha^-(i-1)-\alpha^-(i)\alpha(i)  \; (2\leq i \leq e-1), \; \alpha(e-1)\alpha^-(e-1).
\end{eqnarray*}
\end{lemma}

Theorem \ref{CR}, Theorem \ref{CM} and Lemma \ref{mori} implies the following theorem. 

\begin{theorem} \label{a}
Let $w \in \mathbb{Z}_{\geq 0}$ and $\G$ be any block of $e$-weight $w$ of  $q$-Schur algebra $\mathscr{S}_q (d, d)$. Then the following graded algebra isomorphism holds:  
\begin{equation*}
{\rm HH}^{\ast} (\G) \cong {\rm HH}^{\ast} (A_e^{\otimes w} \rtimes \kk \frks_w ). 
\end{equation*}
\end{theorem}

\begin{remark}
In Theorem \ref{a}, we take a block of a $q$-Schur algebra $\mathscr{S}_q (d, d)$, but we can expand this in the following way: 
\begin{enumerate}
\item[(i)] If $n \geq d$, then $\mathscr{S}_q (n, d)$ is Morita equivalent to $\mathscr{S}_q (d, d)$. 
\item[(ii)] If $n < d$, then we can choose the canonical idempotent (Green's idempotent) $\epsilon \in \mathscr{S}_q (d, d)$ which induces the following isomorphism: 
\begin{equation*}
\mathscr{S}_q (n, d) \cong \epsilon \mathscr{S}_q (d, d) \epsilon.
\end{equation*}
 In this case, $\epsilon$ satisfies the assumption of Corollary \ref{cor1}. Thus we deduce from Corollary \ref{cor1} that the following graded algebra surjection:
\begin{equation*}
{\rm HH}^{\ast}(\mathscr{S}_q (d, d)) \twoheadrightarrow {\rm HH}^{\ast} (\mathscr{S}_q (n, d)).
\end{equation*} 
\end{enumerate}
\end{remark}

\begin{proposition}
Let $(R_{\bullet}, d_{\bullet})$ be the minimal projective resolution of $A_e$ as $A_e$-bimodules. 

\textup{(1)}
We have
\begin{eqnarray*}
&&R_{2s} = \bigoplus_{i = s+1}^e P(i, i) \oplus \left( \bigoplus_{n=1}^s \left( \bigoplus_{j=s-n+1}^{e-2n} P(j, j+2n) \oplus P(j+2n, j)\right)\right),\\
&&R_{2s+1} = \bigoplus_{m=1}^{s+1}\left(\bigoplus_{t=s+2-m}^{e-(2m-1)} P(t, t+2m-1) \oplus P(t+2m-1, t)\right).
\end{eqnarray*}

\textup{(2)}
The differential of $(R_{\bullet},d_{\bullet})$ is given as follows:
\begin{eqnarray*}
d_0 : R_0 &\to& A_e,\\
(i) \otimes (i) &\mapsto& (i),
\end{eqnarray*}
 Let $s \geq0$ and $1 \leq t \leq 2s+1$. With the above notation, for $i \geq 1$, we define the differential $d_i : R_{i} \to R_{i-1}$ recursively as follows:  
\begin{eqnarray*}
d_{4s+1} : R_{4s+1} &\to& R_{4s},\\
(t) \otimes (t+2m-1) &\mapsto& \left\{
\begin{array}{llll}
(-1)^{m+1} \alpha (t-1) \otimes (t+2m-1) \\
+(-1)^m (t) \otimes \alpha^-(t+2m-2) \\
+ (t)\otimes \alpha (t+2m-1) \\
+ \alpha^-(t) \otimes (t+2m-1),  &( 1 \leq m < 2s+1);\\
-(t) \otimes \alpha^-(t+2m-2) \\
+\alpha^-(t) \otimes (t+2m-1),  &(m=2s+1);
\end{array}
\right. \\
(t+2m-1) \otimes (t) &\mapsto& \left\{
\begin{array}{llll}
(-1)^{m+1} (t+2m-1) \otimes \alpha^-(t-1) \\
+(-1)^{m} \alpha(t+2m-2) \otimes (t) \\
+ (t+2m-1)\otimes \otimes \alpha (t) \\
+ \alpha^-(t+2m-1) \otimes (t),  &( 1 \leq m < 2s+1);\\
(t+2m-1) \otimes \alpha(t)\\
-\alpha(t+2m-2) \otimes (t), &(m=2s+1).
\end{array}
\right.
\end{eqnarray*}
Let $s \geq0$ and $1 \leq j \leq 2s+1$. 
\begin{eqnarray*}
d_{4s+2} : R_{4s+2} &\to& R_{4s+1},\\
(i) \otimes (i) &\mapsto& \alpha(i-1) \otimes (i) -(i) \otimes \alpha^-(i-1) \\
&&-(i) \otimes \alpha(i) + \alpha^-(i) \otimes (i), \hspace{4.5em} ( 2s+1 \leq i \leq e);\\
(j) \otimes (j+2n) &\mapsto& \left\{
\begin{array}{lll}
(-1)^n \alpha(j-1) \otimes (j+2n) \\
+ (-1)^{n+1} (j) \otimes \alpha^-(j+2n-1) \\
+\alpha^-(j) \otimes (j+2n) \\
- (j) \otimes \alpha (j+2n), &(1 \leq n < 2s+1);\\
(j) \otimes \alpha^-(j+2n-1)\\
 + \alpha^-(j) \otimes (j+2n), &(n=2s+1); 
\end{array}
\right. \\
(j+2n) \otimes (j) &\mapsto& \left\{
\begin{array}{lll}
(-1)^n \alpha(j+2n-1) \otimes (j) \\
+ (-1)^{n+1} (j+2n) \otimes \alpha^-(j-1) \\
+\alpha^-(j+2n) \otimes (j) \\
- (j+2n) \otimes \alpha(j), &(1 \leq n < 2s+1);\\
- \alpha(J+2n-1) \otimes (j) \\
- (j+2n) \otimes \alpha(j),  &(n=2s+1).
\end{array}
\right. 
\end{eqnarray*}
Let $s \geq0$ and $1 \leq t \leq 2s+2$. 
\begin{eqnarray*}
d_{4s+3} : R_{4s+3} &\to& R_{4s+2},\\
(t) \otimes (t+2m-1) &\mapsto& \left\{
\begin{array}{llll}
(-1)^{m} \alpha (t-1) \otimes (t+2m-1) \\
+(-1)^{m+1} (t) \otimes \alpha^-(t+2m-2) \\
+ (t)\otimes \alpha (t+2m-1) \\
+ \alpha^-(t) \otimes (t+2m-1), &(1 \leq m < 2s+2);\\
-(t) \otimes \alpha^-(t+2m-2)\\
+\alpha^-(t) \otimes (t+2m-1), &(m=2s+2);
\end{array}
\right. \\
(t+2m-1) \otimes (t) &\mapsto& \left\{
\begin{array}{llll}
(-1)^{m} (t+2m-1) \otimes \alpha^-(t-1) \\
+(-1)^{m+1} \alpha(t+2m-2) \otimes (t) \\
+ (t+2m-1)\otimes \alpha (t) \\
+ \alpha^-(t+2m-1) \otimes (t),  &( 1 \leq m < 2s+2);\\
(t+2m-1) \otimes \alpha(t)\\
-\alpha(t+2m-2) \otimes (t), &(m=2s+2).
\end{array}
\right.
\end{eqnarray*}
Let $s \geq0$ and $1 \leq j \leq 2s+2$. 
\begin{eqnarray*}
d_{4s+4} : R_{4s+4} &\to& R_{4s+3},\\
(i) \otimes (i) &\mapsto& -\alpha(i-1) \otimes (i) +(i) \otimes \alpha^-(i-1) \\
&&-(i) \otimes \alpha(i) + \alpha^-(i) \otimes (i), \hspace{3em} ( 2s+1 \leq i \leq e);\\
(j) \otimes (j+2n) &\mapsto& \left\{
\begin{array}{lll}
(-1)^{n+1} \alpha(j-1) \otimes (j+2n) \\
+ (-1)^{n} (j) \otimes \alpha^-(j+2n-1) \\
+\alpha^-(j) \otimes (j+2n) \\
- (j) \otimes \alpha (j+2n), &(1 \leq n < 2s+2);\\
(j) \otimes \alpha^-(j+2n-1) \\
+ \alpha^-(j) \otimes (j+2n), &(n=2s+1); 
\end{array}
\right. \\
(j+2n) \otimes (j) &\mapsto& \left\{
\begin{array}{lll}
(-1)^{n+1} \alpha(j+2n-1) \otimes (j) \\
+ (-1)^{n} (j+2n) \otimes \alpha^-(j-1) \\
+\alpha^-(j+2n) \otimes (j) \\
- (j+2n) \otimes \alpha(j), &(1 \leq n < 2s+2);\\
-(j+2n) \otimes \alpha(j)\\
 - \alpha(j+2n-1) \otimes (j), &(n=2s+2). 
\end{array}
\right. 
\end{eqnarray*}

\end{proposition}

\begin{proof}
We construct the minimal $A_e$-bimodule projective resolution of $A_e$, that is, we construct the following exact sequence:
\begin{equation*}
R_{\bullet} : \cdots \to R_n \xrightarrow{d_n} R_{n-1} \to \cdots \to R_1 \xrightarrow{d_1} R_0 \xrightarrow{d_0} A_e \to 0, 
\end{equation*}
 where $R_n= \bigoplus P(i, j)$ and $P(i, j)$ is the projective $A_e$-bimodule $A_e (i) \otimes (j) A_e$. We obtain from \cite{H} that the projective module $P(i, j)$ occurs in $R_n$ as many times as ${\rm dim}\;{\rm Ext}_{A_e}(S_i, S_j)$. Hence we deduce
\begin{eqnarray*}
&&R_{2s} = \bigoplus_{i = s+1}^e P(i, i) \oplus \left( \bigoplus_{n=1}^s \left( \bigoplus_{j=s-n+1}^{e-2n} P(j, j+2n) \oplus P(j+2n, j)\right)\right),\\
&&R_{2s+1} = \bigoplus_{m=1}^{s+1}\left(\bigoplus_{t=s+2-m}^{e-(2m-1)} P(t, t+2m-1) \oplus P(t+2m-1, t)\right).
\end{eqnarray*}
From this result, we can construct a differential of $A_e$ as above. It is straightforward to check that $d_{\bullet}$ define a complex.

In {\cite{ES}}, they construct an explicit bimodule projective resolution of tame blocks of Hecke algebras of type $A$. Following their method, we show that $d_i$ is the differential of a minimal projective resolution of $A_e$ for any $i \geq 0$. We construct a projective resolution of $A_e/ J A_e$ as left $A_e$-modules, where $J:= J(A_e)$. We denote by $(P_{\bullet}, \delta_{\bullet})$ a projective resolution of $A_e/ J A_e$. Then we have $R_m \otimes_{A_e} A_e/ J A_e \cong P_m$ for all $m \geq 0$ and we deduce that the diagram 
\begin{equation*}
\xymatrix{
\cdots \ar[r] & R_{m+1} \otimes_{A_e} A_e/J A_e  \ar[d]_{\cong} \ar[r]^{d_{m+1} \otimes {\rm id}} & R_{m} \otimes_{A_e} A_e/J A_e  \ar[d]^{\cong} \ar[r]& \cdots &\\
\cdots \ar[r] & P_{m+1} \ar[r]_{\delta_{m+1}} & P_{m} \ar[r] & \cdots. 
 }
\end{equation*}   
commutes for all $m \geq 1$.

We suppose that $\Ker d_{m} \not\subseteq \Im d_{m+1}$ for some $m \geq 1$. Then there exists a non-zero map $\Ker d_m \to \Ker d_m/\Im d_{m+1}$. Therefore there exists a simple $A_e$-bimodule $S \otimes T$ such that $S$ is a left simple $A_e$-module and $T$ is a right $A_e$-module. Moreover there exists a non-zero map $f : \Ker d_{m} \to S \otimes T$. Thus we obtain 
\begin{eqnarray*}
\Im d_{m+1} \otimes_{A_e} A_e / J A_e & \cong & \Im (d_{m+1} \otimes {\rm id} )\\
& \cong & \Ker (d_m \otimes {\rm id} )\\
& \cong & \Ker d_m \otimes_{A_e} A_e/J A_e.  
\end{eqnarray*}
Thus we have
\begin{eqnarray*}
A_e \otimes_{A_e} A_e / J A_e & \xrightarrow{d_{m+1} \otimes {\rm id}} & \Im d_{m+1} \otimes_{A_e} A_e /J A_e \\
& \cong & \Ker d_m \otimes_{A_e} A_e / J A_e \\
& \xrightarrow{f \otimes {\rm id} }& (S \otimes T ) \otimes_{A_e} A_e / J A_e \\
& \cong & S \otimes T.
\end{eqnarray*}
On the other hand, we apply the functor $- \otimes_{A_e} A_e/J A_e$ to the following exact sequence:
\begin{equation*}
R_{m+1} \xrightarrow{d_{m+1}} \Ker d_m \xrightarrow{f} S \otimes T .
\end{equation*}
Then we have
\begin{equation*}
f \circ d_{m+1} = 0.
\end{equation*}
Thus we obtain 
\begin{equation*}
(f \otimes {\rm id}) \circ (d_{m+1} \otimes {\rm id} ) \neq (f \circ d_{m+1} ) \otimes {\rm id}.
\end{equation*}
This is a contradiction. Thus we deduce that $\Ker d_{m} \subseteq \Im d_{m+1}$ for any $m \geq 1$.  
\end{proof}

We give another simple and direct proof of the following Theorem by using the above explicit bimodule projective resolution. 

\begin{theorem}[cf.\ de la Pe{\~n}a-Xi {\cite[Proposition 4.1]{DX}}]  \label{hhgp}
We have the dimension of the Hochschild cohomology group of $A_e$ as follows:
\begin{eqnarray*}
{\rm dim} \; {\rm HH}^n (A_e) = \left\{
\begin{array}{lll}
e, \; & ( n = 0);\\
1, & ( 1 \leq n \leq 2(e-1)); \\
0, & ( 2(e-1) < n).
\end{array}
\right.
\end{eqnarray*}
\end{theorem}

\begin{proof}
If $n=0$, then we have ${\rm dim} \;{\rm HH}^0 (A_e) = {\rm dim} \; Z(A_e) = e$, where $Z(A_e)$ is the center of $A_e$. 

If $1 \leq n \leq \gldim A_e = 2(e-1)$, then there exists a short exact sequence as follows: 
\begin{equation*}
0 \to \Ker d_n \to R_n \to \Ker d_{n-1} \to 0.
\end{equation*}
We apply the functor $\Hom_{\A_}(-, A_e)$ to this exact sequence. Then we have 
\begin{eqnarray*}
0 &\to& \Hom_{\A_} (\Ker d_{n-1} , A_e) \to \Hom_{\A_} (R_n, A_e) \to \Hom_{\A_} (\Ker d_n, A_e) \\
& \to & \Ext^1_{\A_} ( \Ker d_{n-1}, A_e) \to 0.
\end{eqnarray*}
Since $\Ext^i_{\A_} (R_n, A_e)=0$ for all $i \geq 1$, we deduce
\begin{equation*}
\Ext^1_{\A_} (\Ker d_{n-1}, A_e ) \cong \Ext^2_{\A_} ( \Ker d_{n-2}, A_e) \cong \cdots \cong \Ext^{n+1}_{\A_} (A_e, A_e) .
\end{equation*}
We shall determine ${\rm dim}\Hom_{\A_} (R_n, A_e)$ and ${\rm dim}\Hom_{\A_} (\Ker d_n, A_e)$ to determine ${\rm dim}\; {\rm HH}^i (A_e)$. It is straightforward to show ${\rm dim} \Hom_{\A_} (R_n, A_e) = 2e-n-1.$

Next we compute ${\rm dim}\;\Hom_{A_e^{{\rm en}}}(\Ker d_n, A_e)$. We have 
\begin{equation*}
R_{n+1}/ \Ker d_{n+1} \cong \Ker d_n
\end{equation*}
for all $n \geq 0$. Then we have 
\begin{eqnarray*}
\Hom_{\A_} (\Ker d_n, A_e) & \cong& \Hom_{\A_} (R_{n+1}/ \Ker d_{n+1}, A_e)\\
&=& \{ \eta \in \Hom_{\A_}(R_{n+1}, A_e) \; | \; \eta (\Ker d_{n+1} ) = 0. \}.
\end{eqnarray*}
Thus we obtain
\begin{eqnarray*}
{\rm dim} \Hom_{\A_} (\Ker d_n, A_e) = \left\{
\begin{array}{ll}
e-2s-1, \; & ( n= 4s+1, 4s+2);\\
e-2s-2, & ( n= 4s+3, 4s+4 ).
\end{array}
\right.
\end{eqnarray*}
It is clear that ${\rm dim} \;{\rm HH}^n (A_e)=0$ for all $n > \gldim A_e$. 
\end{proof}

\subsection{Main result 2}

From Theorem \ref{hhgp}, we can determine the ring structure of ${\rm HH}^{\ast} (A_e)$.

\begin{theorem}\label{prop}
We have the following $\mathbb{Z}$-graded $\kk$-algebra isomorphism: 
\begin{equation*}
 {\rm HH}^{\ast}(A_e) \cong  \kk [z_{1}, z_{2}, \cdots, z_{e-1}, x, y]/J,
\end{equation*}
where ${\rm deg}\ z_{i}=0$, ${\rm deg}\ x=1$, ${\rm deg}\ y=2$, and 
\begin{eqnarray*}
 J=
\left<{
 \begin{aligned}
z_{i}z_{j},\  z_{i}x, \ z_{i}y, \
x^2, \ xy^{e-1},\ y^e 
\end{aligned}
}
\right>
.
\end{eqnarray*}
\end{theorem}

\begin{proof}
We obtain from Lemma \ref{hhgp} that ${\rm dim}\;{\rm HH}^1 (A_e) = 1$, so we can take  
\begin{equation*}
0 \neq \eta \in {\rm HH}^1 (A_e).
\end{equation*}
Since the Hochschild cohomology ring is super commutative {\cite{G}}, we have 
\begin{equation*}
\eta^2 = 0.
\end{equation*}
Next, we take
\begin{equation*}
0 \neq \theta \in {\rm HH}^2 (A_e).
\end{equation*}
Then, for $1 \leq s \leq e-1$, we show that 
\begin{equation*}
\eta \theta^{s-1}, \theta ^{s} \neq 0
\end{equation*}
by induction on $s$. In the case $s=2$, it is trivial. So we assume that the claim holds for $s>2$. Then we can compute the Yoneda product from an explicit projective resolution of $A_e$, and we have $\theta^{s+1} \neq 0$.
We show that $\eta \theta^{2s} \neq 0$. We can also have an explicit computation of the Yoneda product, and we have $\eta \theta^{2s} \neq 0$. 
\end{proof}

From Proposition \ref{prop}, we obtain the following corollary. 

\begin{corollary}
The even part of ${\rm HH}^{\ast}(A_e)$ is given as follows: 
\begin{equation*}
{\rm HH}^{\rm ev}(A_e)\cong \kk [z_1, \ldots , z_{e-1}, y]/ \langle z_i z_j , z_k y , y^e \; | \; 1 \leq i, j \leq e-1 \rangle , 
\end{equation*}
where  {\rm deg}\ $ z_i =0$, {\rm deg}\ $y =2$.
\end{corollary}

We extend this result to  ${\rm HH}^n(A_e^{\otimes w} \rtimes \kk \mathfrak{S}_{w})$ by using the following result. From the rest of this section, we assume that $l=0$ or $w<l$. 

\begin{proposition} [{Alev-Farinati-Lambre-Solotar \cite[Proposition 3.1]{AFLS}}, {Etingof-Oblomkov \cite[Theorem 3.1]{EO}}] \label{eo}
Let $w \in \mathbb{Z}_{\geq 0}$, and $\Gamma$ be an algebra over $\kk$. Then we have an isomorphism as $\mathbb{Z}_{\geq 0}$-graded vector space over $\kk$. 
\begin{eqnarray*} 
{\rm HH}^{\ast}(\Gamma^{\otimes w} \rtimes \kk \mathfrak{S}_{w}) \cong 
\bigoplus_{ \lambda  \in P_w} \bigotimes_{ i \geq 1} 
({\rm HH}^{\ast}(\Gamma)^{\otimes p_i  (\lambda )})^{\mathfrak{S}_ {p_i  (\lambda )}},
\end{eqnarray*}
where $P_w$ is a set of partitions of $w$ and $p_i (\lambda)$ is a multiplicity of occurrence of $i$ in a partition $\lambda$.
\end{proposition}

In the following, we concentrate on the even part $\kk [z_1, \ldots , z_{e-1}, y]/ \langle z_i z_j , z_k y, y^e \rangle$ of ${\rm HH}^{\ast}(A_e)$. We expand this even part to ${\rm HH}^{{\rm ev}}(A_e^{\otimes w} \rtimes \kk \mathfrak{S}_{w})$ by using Proposition \ref{eo}. We consider the following $\kk$-algebra homomorphism:
\begin{eqnarray*}
\phi : \kk [z_1, \ldots , z_{e-1}, y]/ \langle z_i z_j , z_k y, y^e \rangle &\twoheadrightarrow & \kk [y]/ \langle y^e \rangle \\
 y & \mapsto & y\\
 z_i &\mapsto & 0.
\end{eqnarray*}
Then we have $\Ker \phi = \langle z_i \; |\; 1 \leq i \leq e-1 \rangle$. We consider the $w$-fold tensor product of $\phi$ as commutative ring.  
\begin{equation*}
\phi ^{\otimes w} : ( \kk [z_1, \ldots , z_{e-1}, y]/ \langle z_i z_j , z_k y, y^e \rangle )^{\otimes w} \twoheadrightarrow ( \kk [y]/ \langle y^e \rangle )^{\otimes w}. 
\end{equation*}
We can regard $\phi^{\otimes w}$ as a $\kk \frks_w$-homomorphism, and apply the functor $\Hom_{\kk \frks_w}(\kk, -)$ to $\phi^{\otimes w}$. Then we obtain   
\begin{equation*}
\Hom_{\kk \frks_w}(\kk, \phi^{\otimes w} ) :  ((\kk [z_1, \ldots , z_{e-1}, y]/ \langle z_i z_j , z_k y , y^e \rangle )^{\otimes w} )^{\frks_w} \twoheadrightarrow
(( \kk [y]/ \langle y^e \rangle )^{\otimes w}) ^{\frks_w} .
\end{equation*}
Since
\begin{equation*}
((\kk [y]/ \langle y^e \rangle)^{\otimes w} )^{\frks_w} \cong (\kk [y_1, \ldots , y_w]/ \langle y_1^e, \ldots, y_w^e \rangle )^{\frks_w}, 
\end{equation*}
we consider the kernel of the following natural surjection to give generators of $((\kk [y]/ \langle y^e \rangle)^{\otimes w} )^{\frks_w}$: 
\begin{eqnarray*}
\pi : \Lambda_w  &\to& (\kk [y_1, \ldots , y_w]/ \langle y_1^e, \ldots, y_w^e \rangle )^{\frks_w}\\
x_i & \mapsto & y_i  ,
\end{eqnarray*}
where $\Lambda_w = \kk [x_1, \ldots , x_w]^{\frks_w}$ is the ring of symmetric polynomials in $w$ variables, for more detail on symmetric polynomials see {\cite{M}}.
It is known that
\begin{equation*}
\Ker \pi = \langle p_{e+1}, \ldots, p_{e+w+1} \rangle ,
\end{equation*}
where $p_i$ is the power-sum symmetric polynomial {{\cite[Corollary 3.3]{GA}}}. Consequently, we obtain the following $\mathbb{Z}_{\geq 0}$-graded $\kk$-algebra isomorphism:
\[
\mathrm{HH}^{{\rm ev}}\left(A_{e}^{\otimes w}\rtimes\mathbf{k}\mathfrak{S}_{w}\right) / \Ker \Hom_{\kk \mathfrak{S}_w}(\kk, \phi^{\otimes w})\simeq\Lambda_{w}/\left\langle p_{e+1},\cdots,p_{e+w+1}\right\rangle.
\]

\section*{acknowledgement}
The author is greatly indebted to Hyohe Miyachi, Akihiro Tsuchiya and Yoshiyuki Kimura for their helpful advices. This work is supported by by Grant-in-Aid for JSPS Fellowships No. H15J09492.

\bibliographystyle{alpha}
\bibliography{Bib_File1}

\end{document}